\documentclass[12pt]{article}
\usepackage{amsmath}
\usepackage{latexsym}
\usepackage{amssymb}
\newtheorem{thm}{Theorem}[section]
\newtheorem{la}[thm]{Lemma}
\newtheorem{Defn}[thm]{Definition}
\newtheorem{Remark}[thm]{Remark}
\newtheorem{Note}[thm]{Note}
\newtheorem{prop}[thm]{Proposition}

\newtheorem{cor}[thm]{Corollary}
\newtheorem{Example}[thm]{Example}
\newtheorem{Examples}[thm]{Examples}
\newtheorem{Problems}[thm]{Problems}

\newtheorem{Problem}[thm]{Problem}
\newtheorem{Number}[thm]{\!\!}
\newenvironment{defn}{\begin{Defn}\rm}{\end{Defn}}

\newenvironment{example}{\begin{Example}\rm}{\end{Example}}

\newenvironment{rem}{\begin{Remark}\rm}{\end{Remark}}
\newenvironment{numba}{\begin{Number}\rm}{\end{Number}}
\newenvironment{proof}{{\noindent\bf Proof.}}%
                  {\nopagebreak\hspace*{\fill}$\Box$\medskip\medskip\par}   
\newcommand{\Punkt}{\nopagebreak\hspace*{\fill}$\Box$}

\newcommand{\wb}{\overline}

\DeclareMathOperator{\tr}{tr}

\newcommand{\tensor}{\otimes}
\newcommand{\impl}{\Rightarrow}

\newcommand{\mto}{\mapsto}

\newcommand{\Ad}{\mbox{\rm A\hspace*{-.2 mm}d}}
\newcommand{\N}{{\mathbb N}}

\newcommand{\bO}{{\mathbb O}}
\newcommand{\F}{{\mathbb F}}
\newcommand{\K}{{\mathbb K}}

\newcommand{\Q}{{\mathbb Q}}

\newcommand{\Z}{{\mathbb Z}}

\newcommand{\cg}{{\mathfrak g}}

\newcommand{\co}{{\mathfrak o}}

\DeclareMathOperator{\Aut}{Aut}

\newcommand{\one}{{\bf 1}}
\newcommand{\sub}{\subseteq}
\DeclareMathOperator{\GL}{GL}
\DeclareMathOperator{\Orth}{O}
\DeclareMathOperator{\SL}{SL}
\DeclareMathOperator{\UT}{UT}

\DeclareMathOperator{\fin}{fin}
\newcommand{\gl}{\mbox{${\mathfrak g}{\mathfrak l}$}}
\newcommand{\ut}{\mbox{${\mathfrak u}{\mathfrak t}$}}
\newcommand{\Sl}{\mbox{${\mathfrak s}{\mathfrak l}$}}

\DeclareMathOperator{\id}{id}

\newcommand{\cB}{{\mathcal B}}

\newcommand{\cH}{{\mathcal H}}

\DeclareMathOperator{\Sym}{Sym}
\DeclareMathOperator{\Inn}{Inn}

\DeclareMathOperator{\dt}{det}
\DeclareMathOperator{\car}{char}

\newcommand{\dotarrow}{\mbox{$\,\;\cdots\hspace*{-.63mm}\gtrdot\;\,$}}
\begin{document}
\renewcommand{\thefootnote}{\fnsymbol{footnote}}
\begin{center}
{\Large\bf
Directions of automorphisms of Lie groups\vspace{2mm}
over local fields
compared to\vspace{2mm}
the directions of Lie algebra automorphisms}\\[5mm]
{\bf Helge Gl\"{o}ckner and George A. Willis\footnote{Research
supported by DFG grant 447 AUS-113/22/0-1
and ARC grant LX 0349209.}}\vspace{2mm}
\end{center}
\renewcommand{\thefootnote}{\arabic{footnote}}
\setcounter{footnote}{0}
\begin{abstract}
\hspace*{-7.2 mm}
To each totally disconnected,
locally compact topological group~$G$
and each group~$A$ of automorphisms
of~$G$,
a pseudo-metric space $\partial A$
of ``directions''
has been associated by U. Baumgartner
and the second author.
Given a Lie group $G$ over a local field,
it is a natural idea to try to define a map
\[
\Phi\colon
\partial \Aut_{C^\omega}(G)\to \partial \Aut(L(G))\,,
\quad \partial \alpha\mto \partial L(\alpha)
\]
which takes the direction of an analytic automorphism
of~$G$ to the direction of the associated Lie algebra
automorphism.
We show that, in general,
$\Phi$ is not well-defined.
Also, it may happen that
$\partial L(\alpha)=\partial L(\beta)$
although $\partial \alpha\not=\partial\beta$.
However, such pathologies are absent for
a large class of groups:
we show that
$\Phi\colon \partial\hspace*{.1mm}\Inn(G)\to \partial\hspace*{-.2mm}\Aut(L(G))$
is a well-defined isometric embedding
for each generalized Cayley group~$G$.
Some counterexamples concerning the existence
of small joint tidy subgroups
for flat groups of automorphisms are also provided.\vspace{2.5mm}

\end{abstract}
{\footnotesize {\em Classification}:
22D05 (primary); 
20G25; 
22D45; 
22E15; 
22E35\\[1mm] 
{\em Key words}: totally disconnected group; automorphism;
direction; local field; Lie group; algebraic group;
Cayley transform; Cayley group; generalized Cayley group;
scale; scale function;
tidy subgroup; flat group; tidy automorphism;
small subgroup}\vspace{2.5mm}
\begin{center}
{\large\bf Introduction}
\end{center}
In a recent article~\cite{BaW},
U. Baumgartner and the second author
associated a pseudo-metric space $\partial A$
of ``directions'' to
each totally disconnected, locally compact group~$G$
and group~$A$ of (bicontinuous) automorphisms of~$G$.
The completion $\overline{\partial\hspace*{.15mm} G}$
of the metric
\vspace{1.6mm}space associated with
the space
$\partial \hspace*{.2mm}G:=\partial \Inn(G)$
of directions of inner automorphisms
generalizes familiar objects.
For example,
$\overline{\partial \hspace*{.2mm}G}$ is homeomorphic to the
spherical building at infinity
if~$G$ is a semisimple group over a local field~\cite{BaW}.
We recall from~\cite{BaW}
that an automorphism $\alpha$ of~$G$
is said to \emph{move to infinity}
if, for all
compact open subgroups $V\subseteq W$
of $G$, there exists a positive integer~$n$
such that $\alpha^n(V)\not\subseteq W$.
Each automorphism $\alpha\in A$
which moves to infinity
can be assigned a ``direction''
$\partial\alpha\in \partial A$,
and every element
of $\partial A$ arises
in this way.
The first part of the article
is devoted to the relations between
the space of directions
$\partial\hspace*{-.35mm}\Aut_{C^\omega}(G)$
of analytic automorphisms
of a Lie group~$G$ over a local field~$\K$
and the space of directions
$\partial\hspace*{-.35mm}\Aut(L(G))$ of automorphisms
of (the additive group of)
its Lie algebra.
It is a natural idea to try to define
a map
\[
\Phi\colon
\partial \hspace*{-.2mm}\Aut_{C^\omega}(G)
\to \partial \hspace*{-.2mm}\Aut(L(G))\,,
\quad \partial \alpha\mto \partial L(\alpha)
\]
which takes the direction of an analytic automorphism
of~$G$ to the direction of the associated linear
automorphism.
Section~\ref{seccount}
compiles negative results:
We show that, in general,
$\Phi$ is not well-defined,\footnote{Not even as a map into
$\wb{\partial\hspace*{-.3mm}\Aut(L(G))}$.}
at least if \mbox{$\car(\K)>0$}\linebreak
(Example~\ref{ex2}).
Furthermore, it may happen that
$\partial L(\alpha)=\partial L(\beta)$
although $\partial \alpha\not=\partial\beta$
(see Examples~\ref{ex1} and~\ref{ex3}).
Section~\ref{secposi} is devoted
to positive results:
we describe additional conditions
ensuring that $\Phi\colon \partial G\to \partial \!\Aut(L(G))$
is well-defined, respectively, an isometric embedding
(Proposition~\ref{extracond}).
These conditions are satisfied by
a large class of linear algebraic groups
(called\linebreak
 ``generalized Cayley groups'' here),
for which a well-behaved
analogue of the Cayley transform
is available.
These groups are variants of the Cayley groups
investigated recently in~\cite{LPR}.
In particular, we shall see in Corollary~\ref{orthog}
that $\Phi\colon \partial G\to \partial \Aut(L(G))$
is a well-defined isometric embedding if~$G$
is a general linear group,
a special linear group (if $\car(\K)=0$),
an orthogonal group (if~$\car(\K)\not=2$)
or the group of all invertible upper
triangular $n\!\times\! n$-matrices.\\[2.5mm]
We recall two
concepts from the second author's
structure theory of totally disconnected
groups (see \cite{Wi1}--\cite{Wi2b}):
Given $\alpha\in \Aut(G)$,
its \emph{scale} is defined as
the minimum index $s_G(\alpha):=
\min_U \,[U: U\cap\alpha^{-1}(U)] \in \N$,
where~$U$ ranges through the set
of compact open subgroups
of~$G$.
If the minimum is attained at~$U$,
then the compact open subgroup~$U$
is called \emph{tidy for~$\alpha$}
(see~\cite{Wi2}; cf.\ \cite{Wi1}
for equivalent earlier definitions).\\[2.5mm]
These concepts play an important role
in the construction of the space
of directions.
Notably, the scale is needed
to define the pseudo-metric~on~$\partial\hspace*{-.3mm}\Aut(G)$.
The scale also facilitates
a convenient characterization:
An automorphism
$\alpha\in \Aut(G)$ moves to infinity if and only if
$s_G(\alpha)>1$ (see \cite[Lemma~3]{BaW}).\\[2.5mm]
Following~\cite{TID},
an automorphism~$\alpha$ of~$G$
is called \emph{tidy}
if $G$ has arbitrarily small
compact open subgroups which are tidy for~$\alpha$.
Tidiness of automorphisms is a useful regularity property,
which rules out many pathologies
(see~\cite{TID}, \cite{POS}).
Also for our present purposes,
tidiness is of interest:
it ensures that $\alpha$ moves to infinity
if and only if so does~$L(\alpha)$
(see Lemma~\ref{liescale}\,(c) and Lemma~\ref{toinfty}\,(b)).\\[2.5mm]
We recall from~\cite{Wi2b}
that a group $\cH\leq \Aut(G)$
of automorphisms is called
\emph{flat} if~$G$ has a compact, open subgroup
which is tidy for each $\alpha\in \cH$.
For some purposes,
flat groups in totally disconnected
groups can be used as substitutes
for tori in algebraic groups.
They are also used
in the proof that the space of directions
of a semisimple algebraic group~$G$
over a local field
is homeomorphic to the spherical building
of~$G$ (see \cite{BaW}),
as mentioned above.\\[2.5mm]
The second part of this article
(Section~\ref{secsmall}) provides
new results concerning flat groups
of automorphisms
of a totally disconnected group~$G$.
Notably, we describe examples
of flat groups $\cH$ of tidy automorphisms
such that $G$ does not have small subgroups
which are tidy for all $\alpha\in \cH$
simultaneously.
\section{Notation and auxiliary results}\label{secprel}
In this section, we recall some definitions
and basic facts
which are necessary to define and discuss
the space of directions $\partial \Aut(G)$
of a totally disconnected, locally compact group~$G$.
We also compile some facts concerning the scale
of automorphisms of Lie groups
over local fields.
\begin{numba}
On the set of compact open subgroups
of a totally disconnected, locally compact group~$G$,
a metric can be defined using the formula
$d(V,W):=d_+(V,W)+d_+(W,V)$,
where $d_+(V,W):= \log[V:V\cap W]$ (see \cite[p.\,395]{BaW}).
Two automorphisms~$\alpha$ and~$\beta$
of $G$ are called \emph{asymptotic} if the sequence
$(d(\alpha^{nk}(V),\beta^{n\ell}(W)))_{n\in {\mathbb N}}$
is bounded for some (hence any)
compact open subgroups $V,W\subseteq G$
and certain $k,\ell\in {\mathbb N}$
(cf.\ Definition~6 and Lemma~7 in~\cite{BaW}).
In this case, we write $\alpha\asymp \beta$.
Given a group $A\leq \Aut(G)$ of automorphisms,
being asymptotic is an equivalence relation
on the set~$A_>$ of automorphisms moving
towards infinity~\cite[Lemma~10]{BaW}.
The equivalence class of $\alpha\in A_>$ is
called its \emph{direction}.
We write $\partial A:=\{\partial \alpha\colon \alpha\in A_>\}$
and abbreviate $\partial G:=\partial \hspace*{-.3mm}\Inn(G)$,
where $\Inn(G)$ is the group of all
inner automorphisms of~$G$.
To construct a pseudo-metric on~$\partial A$,
as an auxiliary notion define
\[
\delta^{V,W}_n(\alpha,\beta) := 
\min\left\{
\frac{d_+(\alpha^n(V),\beta^k(W))}{n\log(s_G(\alpha))}\colon
\mbox{$k\in \N$ such that $s_G(\beta)^k\leq s_G(\alpha)^n$}
\right\}
\]
for $\alpha,\beta\in A_>$,
compact open subgroups $V,W\sub G$
and $n\in \N$. It can be shown that
\[
\delta_+(\alpha,\beta)\; :=\; \limsup_{n\to\infty}
\delta_n^{V,W}(\alpha,\beta)\;\in\; [0,1]\,,
\]
and that $\delta_+(\alpha,\beta)$ is independent of the
choice of~$V$ and~$W$ (see \cite[p.\,406]{BaW}).
Furthermore,
\[
\delta(\alpha,\beta)\; :=\;
\delta_+(\alpha,\beta)+\delta_+(\beta,\alpha)
\]
defines a pseudo-metric~$\delta$ on~$A_>$
by \cite[Corollary~16]{BaW},
and a well-defined
pseudo-metric on $\partial A$ (also denoted~$\delta$)
is obtained via $\delta(\partial \alpha,\partial\beta):=
\delta(\alpha,\beta)$ (see \cite[Lemma~17]{BaW}).
The completion of the metric space $\partial A/\delta^{-1}(0)$
associated with~$\partial A$ is denoted by
$\wb{\partial A}$.
Finally, we let $\Delta_G(\alpha)$
be the module of $\alpha\in\Aut(G)$,
defined as
$\Delta_G(\alpha):=\frac{\mu(\alpha(U))}{\mu(U)}$,
where~$\mu$ is a Haar measure on~$G$ and~$U\sub G$
a non-empty, relatively compact, open set.
\end{numba}
\begin{numba}
All Lie groups considered in this article
are finite-dimensional.
We use ``$C^\omega$''
as a shorthand for ``analytic.'' 
Beyond $C^\omega$-Lie groups
over a local field~$\K$
(as in~\cite{Ser}),
it is possible to define
$C^k$-Lie groups over~$\K$
for each $k\in \N\cup\{\infty\}$,
based on a notion of
$C^k$-map between open subsets
of finite-dimensional
(or topological)
vector spaces introduced
in~\cite{BGN}.\footnote{By
\cite[Lemmas~3.2 and~4.4]{IMP}, a map is~$C^1$ if and only
if it is strictly differentiable
at each point in the sense of~\cite[1.2.2]{BVA}.
For $C^k$-maps on subsets of~$\K$,
see already~\cite{Sch}.}
Every $C^\omega$-Lie group is also a $C^k$-Lie group,
but the converse is valid
only in zero characteristic (see \cite{NAN}, \cite{ANA}).
For the sake of added generality,
we therefore formulate
our results for $C^1$-automorphisms
of $C^1$-Lie groups.
However, all of our concrete examples
will be $C^\omega$-Lie groups,
and readers unfamiliar with
$C^k$-maps over local fields
are invited to
replace ``$C^1$''
by ``$C^\omega$''
throughout the article.
\end{numba}
In the next two lemmas, $\K$ is
a local field, $\wb{\K}$ an algebraic closure of~$\K$
and $|.|\colon \wb{\K}\to [0,\infty[$
an absolute value whose restriction
to $\K^\times$ is $\text{mod}_\K$
(cf.\ \cite[Chapter~I]{Wei} and \cite[\S\,14]{Sch}).
``Tidiness'' and ``$s_G$''
are as in the Introduction.
\begin{la}\label{liescale}
Let $\alpha$ be a $C^1$-automorphism of a $C^1$-Lie group~$G$
over a local field~$\K$,
and $L(\alpha):=T_1(\alpha)$
be the corresponding linear automorphism
of the tangent space $L(G):=T_1(G)$.
Then the following holds:\vspace{.5mm}

\begin{itemize}
\item[\rm(a)]
Let
$\lambda_1,\ldots,\lambda_{\,\dim \!L(G)}\in \wb{\K}$
be the eigenvalues of $\alpha\tensor_\K\wb{\K}\in
\GL(L(G)\tensor_\K \wb{\K})$,
with repetitions
according to the algebraic multiplicities.
Then
\begin{equation}\label{formsca}
s_{L(G)}(L(\alpha))\;\; =\; \prod_{j\;
\mbox{\rm\footnotesize such that}\;|\lambda_j|>1}
|\lambda_j|\,.
\end{equation}
\item[\rm(b)]
$s_G(\alpha)$ divides $s_{L(G)}(L(\alpha))$,
whence
$s_G(\alpha)\leq s_{L(G)}(L(\alpha))$ in particular.
\item[\rm(c)]
$s_G(\alpha)= s_{L(G)}(L(\alpha))$
holds if and only if~$\alpha$ is
tidy.
\item[\rm (d)]
If $\alpha$ is an inner automorphism
and there exists an injective, continuous
homomorphism $f\colon G\to \GL_n(\K)$ for some~$n$,
then
$s_G(\alpha)=s_{L(G)}(L(\alpha))$.
\item[\rm(e)]
If $\car(\K)=0$ or if
$G$ is one-dimensional, then $s_G(\alpha)=s_{L(G)}(L(\alpha))$.
\end{itemize}
\end{la}
\begin{proof}
See \cite{POS} for (a)--(d)
and the first
half of~(e). To prove the second
half of~(e), assume that $\dim_\K(L(G))=1$.
If $s_{L(G)}(L(\alpha))=1$, then also
$s_G(\alpha)=1$, by~(b).
Now assume that
$s_{L(G)}(L(\alpha))>1$.
Since $L(G)$ is $1$-dimensional,
we have $L(\alpha)=\lambda\,\id_{L(G)}$
for some $\lambda\in \K^\times$.
Then $|\lambda|=s_{L(G)}(L(\alpha))>1$,
by~(\ref{formsca}).\hspace*{-.3mm}
Hence $L(\alpha^{-1})\! =\! \lambda^{-1}\id_{L(G)}$
with $|\lambda^{-1}|\! <\! 1$ and thus
\mbox{$s_{L(G)}(L(\alpha^{-1}))\!=\!1$,}
entailing that also
$s_G(\alpha^{-1})=1$.
Therefore
\[
s_G(\alpha)=
\frac{s_G(\alpha)}{s_G(\alpha^{-1})}=\Delta_G(\alpha)
=\Delta_{L(G)}(L(\alpha))
=\frac{s_{L(G)}(L(\alpha))}{s_{L(G)}(L(\alpha^{-1}))}=s_{L(G)}(L(\alpha))\, ,
\]
using
\cite[p.\,354, Corollary~1]{Wi1}
and \cite[Chapter~III, \S\,3.16, Proposition~55]{Bo2}.\vspace{-1mm}
\end{proof}
We remark that if
$\car(\K)=0$ or
$G$ is a Zariski-connected reductive
algebraic group over~$\K$
and~$\alpha$ an inner automorphism,
then $s_G(\alpha)$ has been
expressed by the right hand side
of~(\ref{formsca}) already in~\cite{SCA}
and \cite{CON}, respectively.\\[1.3mm]
Since $\alpha$ moves to infinity
if and only if $s_G(\alpha)>1$,
Lemma~\ref{liescale} implies:
\begin{la}\label{toinfty}
Let $\alpha$ be a $C^1$-automorphism of a $C^1$-Lie group~$G$
over a local field~$\K$.
Then the following holds:
\begin{itemize}
\item[\rm (a)]
If $\alpha$ moves to infinity, then $L(\alpha)$ moves to infinity.
\item[\rm (b)]
If $s_G(\alpha)=s_{L(G)}(L(\alpha))$,
then $\alpha$ moves to infinity if and only
if $L(\alpha)$ moves to infinity.
\item[\rm (c)]
$L(\alpha)$ moves to infinity if and only if
$L(\alpha)$ has an eigenvalue $\lambda$ in $\wb{\K}$
of absolute value $|\lambda|>1$.
\end{itemize}
\end{la}
\begin{proof}
Part~(a) is an immediate consequence of Lemma~\ref{liescale}\,(b).
Part~(b) is obvious, and (c) follows from
Lemma~\ref{liescale}\,(a).
\end{proof}
\section{Counterexamples concerning directions}\label{seccount}
In this section, we provide examples of Lie groups
over local fields with pathological properties,
as announced in the introduction.\\[3mm]
The following simple observation is useful
for our discussions.
\begin{la}\label{trivobs}
Let $G$ be a totally disconnected,
locally compact group,
$V$ and~$W$ be compact, open subgroups of~$G$
and
$\alpha,\beta\in \Aut(G)_>$.
If $s_G(\alpha)=s_G(\beta)$
and $\beta(W)\supseteq W$, then
\[
\delta_n^{V,W}(\alpha,\beta)\;=\;
\frac{d_+(\alpha^n(V),\beta^n(W))}{n\log(s_G(\alpha))}
\quad\mbox{for each $n\in \N$.}
\]
\end{la}
\begin{proof}
This is clear from the definitions.\vspace{-5mm}
\end{proof}
\begin{numba}
In our first and second example,
we consider Lie groups over the field
$\K:=\F(\!(X)\!)$ of formal Laurent series over a finite field~$\F$
of~$q$ elements. The corresponding ring of formal
power series will be abbreviated $\bO:=\F[\![X]\!]$.
The following notations are useful:
Given $z=\sum_{k=-\infty}^\infty a_kX^k\in \K$,
where $(a_k)_{k\in \Z}\in \F^{(-\N)}\times \F^{\N_0}$,
we define
\[
\;z^{(1)} :=  \sum_{k=1}^\infty a_kX^k\, ; \quad\;\;\,
z^{(2)} :=  \sum_{k=1}^\infty a_{-k}X^{-k}\, ;\qquad\qquad\;\,
z^{(3)} := a_0 X^0\, ;\qquad\quad\;
\]
\[
z^{(4)} :=  \sum_{k=0}^\infty a_k X^k ;\quad\;\;\;
z^{(5)} := \sum_{k=0}^\infty a_{-(2k+1)}X^{-(2k+1)};\quad
z^{(6)} := \sum_{k=1}^\infty a_{-2k}X^{-2k}.
\]
Thus $z=z^{(1)}+z^{(2)}+z^{(3)}=z^{(4)}+z^{(5)}+z^{(6)}$.
\end{numba}
\begin{example}\label{ex2}
We describe $\K$-analytic automorphisms $\alpha,\beta\in \Aut(G)_>$ of
the Lie group $G:=(\K^2, +)$
such that $\partial\alpha=\partial \beta$ but
$\delta(\partial L(\alpha),\partial L(\beta))>0$,
whence $\partial L(\alpha)\not=\partial L(\beta)$
in particular.\\[3mm]
Let $\alpha$ be the (linear) automorphism
$\alpha\colon G\to G$, $\alpha(v,w):=(X^{-1}v,X^{-1}w)$ for
$v,w\in \K$, and define
$\beta\colon G\to G$ via
\[
\beta(v,w)\, :=\,
\big( v^{(1)} + X^{-1}(v^{(2)}+v^{(3)}) + w^{(3)} ,\,
X^{-2}w^{(1)}+ X^{-1}w^{(2)}\big)\,.
\]
It is clear that $\beta$ is bijective and a homomorphism of groups;
since $\beta$ coincides with the linear isomorphism $\K^2\to\K^2$,
$(v,w)\mto (v,X^{-2}w)$
on the open zero-neighbourhood $(X \bO)\times (X \bO)$,
we deduce that $\beta$ is a $\K$-analytic automorphism of~$G$
and $L(\beta)(v,w)=(v,X^{-2}w)$.\\[3mm]
To see that
$\delta(\partial L(\alpha),\partial L(\beta))>0$,
note first that
$s_{L(G)}(L(\alpha))=|X^{-1}|^2>0$
and $s_{L(G)}(L(\beta))=|X^{-2}|>1$
by Lemma~\ref{liescale}\,(a),
whence both~$L(\alpha)$ and~$L(\beta)$ move to infinity.
For each $n\in \N$, we have
$L(\alpha)^n(\bO\times\bO)=(X^{-n}\bO)\times (X^{-n}\bO)$
and $L(\beta)^n(\bO\times\bO)=\bO\times X^{-2n}\bO$
with intersection $\bO\times X^{-n}\bO$, whence
\begin{equation}\label{lle}
d_+(L(\alpha)^n(\bO^2),L(\beta)^n(\bO^2))\; = \;
\log \,[X^{-n}\bO:\bO]\; =\; \log(q^n)\;=\; n\log(q)\,.
\end{equation}
Since $s_{L(G)}(L(\alpha))=s_{L(G)}(L(\beta))$
and $L(\beta)(\bO\times\bO)=\bO\times X^{-2}\bO\supseteq \bO\times\bO$,
Lemma~\ref{trivobs} can be applied. Combined with~(\ref{lle}),
it shows that
\[
\delta^{\bO^2,\bO^2}_n(L(\alpha),L(\beta))\;=\;
\frac{d_+(L(\alpha)^n(\bO^2),L(\beta)^n(\bO^2))}{n\log s_{L(G)}(L(\alpha))}
\; =\;
\frac{\log(q)}{\log s_{L(G)}(L(\alpha))}.
\]
Letting
$n\to\infty$, we infer that
$\delta_+(L(\alpha),L(\beta))=\frac{\log(q)}{\log s_{L(G)}(L(\alpha))}>0$.
Hence $\delta(\partial L(\alpha),\partial L(\beta))
=\delta_+(L(\alpha),L(\beta))
+\delta_+(L(\beta),L(\alpha))>0$.\\[3mm]
To see that $\alpha$ moves to infinity, let $V\sub W$
be compact open subgroups
of~$G$. There exist $k,\ell\in \N_0$
such that $X^k\bO\times X^k\bO\sub V$
and $W\sub X^{-\ell}\bO\times X^{-\ell}\bO$.
Set $n:=k+\ell+1$.
Then $v:=(0,X^k)\in V$ but
$\alpha^n(v)=(0,X^{-\ell-1})\not\in W$,
whence $\alpha^n(V)\not\sub W$.
An analogous argument shows that
$\beta$ moves to infinity.
To complete our discussion,
note that $\alpha^n(\bO\times\bO )=(X^{-n}\bO)\times (X^{-n}\bO)
=\beta^n(\bO\times\bO)$
for all $n\in \N$, whence $d(\alpha^n(\bO^2),\beta^n(\bO^2))=0$
for all $n\in \N$ and thus $\partial\alpha=\partial\beta$.
\end{example}
\begin{example}\label{ex1}
We describe
$\K$-analytic automorphisms $\alpha, \beta\in \Aut(G)_>$
of the one-dimensional Lie group $G:=(\K,+)$
such that
$\delta(\partial \alpha, \partial \beta)>0$
(and hence $\partial \alpha\not=\partial \beta$),
but
$L(\alpha)=L(\beta)$ and thus $\partial L(\alpha)=\partial L(\beta)$.\\[3mm]
Consider the linear automorphism
$\alpha\colon G\to G$, $\alpha(z):=X^{-1}z$
and the~map
\[
\beta \colon G\to G\, ,
\quad \beta(z)\, := \, X^{-1}z^{(4)}+ X^{-2}z^{(5)}+ z^{(6)}\,.
\]
Then~$\beta$ is bijective and a homomorphism,
and from $\beta|_\bO=\alpha|_\bO$ we
now deduce that $\beta$ is a $\K$-analytic automorphism with
$L(\beta)=L(\alpha)$.
By Lemma~\ref{toinfty}\,(c),
$L(\beta)=L(\alpha)\colon z\mto X^{-1}z$
moves to infinity
and hence so do~$\alpha$
and~$\beta$, by Lemma~\ref{liescale}\,(e)
and Lemma~\ref{toinfty}\,(b).
To see that $\delta(\partial \alpha,\partial \beta)>0$,
note that $\alpha^n(\bO)=X^{-n}\bO$ and
$\beta^n(\bO)=\sum_{k=0}^{n-1}\F X^{-(2k+1)}+\bO$
for $n\in \N$.
For $n$ odd, say $n=2\ell+1$,
this entails that
\[
\alpha^n(\bO)\cap \beta^n(\bO)\;=\;
\sum_{k=0}^\ell\F X^{-(2k+1)}+\bO\,,
\]
whence
$[\alpha^n(\bO): \alpha^n(\bO)\cap \beta^n(\bO)] = q^\ell$
and hence
\begin{equation}\label{pree}
d_+(\alpha^n(\bO),\beta^n(\bO))
\;=\; \log(q^\ell)\; =\; \ell\,\log(q)\,.
\end{equation}
Since $s_G(\alpha)=s_{L(G)}(L(\alpha))=
s_{L(G)}(L(\beta))=s_G(\beta)$
by Lemma~\ref{liescale}\,(e)
and $\beta(\bO)= \F X^{-1}+\bO\supseteq \bO$,
we can apply Lemma~\ref{trivobs} to see that
\begin{eqnarray*}
\delta_+(\alpha,\beta)
& = &
\limsup_{n\to\infty}
\frac{d_+(\alpha^n(\bO),\beta^n(\bO))}{n\log(s_G(\alpha))}
\;\,\geq\;\,
\limsup_{\ell\to\infty}
\frac{d_+(\alpha^{2\ell+1}(\bO),
\beta^{2\ell+1}(\bO))}{(2\ell+1)\log(s_G(\alpha))}\\
& = &
\lim_{\ell\to\infty}\frac{\ell\, \log(q)}{(2\ell+1)\log(s_G(\alpha))}
\;\, =\;\, \frac{\log(q)}{2 \log(s_G(\alpha))}\;>\; 0\,,
\end{eqnarray*}
where we used~(\ref{pree}) to pass to the second line.
Hence $\delta(\partial\alpha,\partial\beta)>0$.
\end{example}
\begin{example}\label{ex3}
We now describe a one-dimensional $p$-adic Lie group~$G$
and analytic automorphisms $\alpha, \beta\in \Aut(G)_>$
such that $\delta(\partial \alpha, \partial \beta)>0$
(and hence $\partial \alpha\not=\partial \beta$),
but $L(\alpha)=L(\beta)$ and thus
$\partial L(\alpha)=\partial L(\beta)$.\\[3mm]
A suitable $p$-adic Lie group is
$G:=\Q_p\times (\Q_p/\Z_p)$;
the desired automorphisms are
$\alpha\colon G\to G$, $\alpha(x,y):=(p^{-1}x,y)$
and
\[
\beta\colon G\to G\, ,\quad \beta(x,y)\,
:=\, \big( p^{-1}x,\, y+q(p^{-1}x) \big)\,,
\]
where $q\colon \Q_p\to \Q_p/\Z_p$ is the quotient homomorphism.
It is clear that $\alpha$ and $\beta$ are
automorphisms; their inverses are given
by $\alpha^{-1}(x,y)=(px,y)$ and $\beta^{-1}(x,y)=(px,y-q(x))$,
respectively.
Since~$\alpha$ and~$\beta$ coincide on the open zero-neighbourhood
$p \Z_p \times\{0\}$, we have
$L(\beta)=L(\alpha)\colon \Q_p\to\Q_p$, $x\mto p^{-1}x$.
Using Lemma~\ref{toinfty}\,(c), we deduce that
$L(\alpha)=L(\beta)$ moves
to infinity.
By Lemma~\ref{liescale}\,(e)
and Lemma~\ref{toinfty}\,(b),
also~$\alpha$ and~$\beta$ move to infinity.\\[3mm]
We claim that
$\alpha^n(V)\cap \beta^n(V)=V$
for each $n\in \N$,
where $V:=\Z_p \times\{0\}$.
If this is true, then
$[\alpha^n(V): \alpha^n(V)\cap\beta^n(V)]
= [p^{-n}\Z_p \times \{0\}: \Z_p \times\{0\}]
= [p^{-n}\Z_p : \Z_p ] = p^n$
and thus
\[
d_+(\alpha^n(V),\beta^n(V))\; =\; n\log(p)\, .
\]
Since $s_G(\alpha)=s_{L(G)}(L(\alpha))=s_{L(G)}(L(\beta))=
s_G(\beta)$ by Lemma~\ref{liescale}\,(e)
and $\beta(V)\supseteq
\beta(p\Z_p \times\{0\})=\Z_p \times\{0\}=V$,
we can apply Lemma~\ref{trivobs}
to see that
\[
\delta_n^{V,V}(\alpha,\beta)\;=\;
\frac{d_+(\alpha^n(V),\beta^n(V))}{n\log(s_G(\alpha))}\;=\;
\frac{\log(p)}{\log (s_G(\alpha))}
\]
for each $n\in\N$.
Hence $\delta(\partial\alpha,\partial\beta)\geq
\delta_+(\alpha,\beta)=\frac{\log(p)}{\log (s_G(\alpha))}>0$.
It only remains to prove the claim.
If $x=\sum_{k=0}^\infty a_k p^k\in\Z_p$
with $a_k\in \{0,1,\ldots, p-1\}$,
then
\begin{equation}\label{explict}
\beta^n(x,0)\,=\,
\Big( p^{-n}x,\,
\sum_{k=1}^n\big( \sum_{j=0}^{n-k}
a_j\big) p^{-k}\, + \Z_p\Big)
\end{equation}
for each $n\in \N$, by a simple induction.
If this element is in $p^{-n}\Z_p \times \{0\}$, then
$\sum_{k=1}^n\big( \sum_{j=0}^{n-k}
a_j\big) p^{-k} \in \Z_p$.
Hence
$a_0p^{-n}\equiv \sum_{k=1}^n\big( \sum_{j=0}^{n-k}
a_j\big) p^{-k} \equiv 0$ modulo $p^{-(n-1)}\Z_p$
and thus $a_0=0$.
Therefore
$\sum_{k=1}^{n-1}\big( \sum_{j=1}^{n-k}
a_j\big) p^{-k} \in \Z_p$.
Repeating the argument, we find that
$a_0=a_1=\cdots=a_{n-1}=0$ and thus 
\mbox{$\beta^n(x,0)\in \Z_p \times \{0\}$.}
We have shown that
$\alpha^n(\Z_p \times\{0\})\cap \beta^n(\Z_p \times\{0\})
\sub \Z_p \times\{0\}$.
On the other hand, $\beta^n(\Z_p \times \{0\})
\supseteq \beta^n((p^n\Z_p)\times\{0\})=\Z_p\times\{0\}$.
Hence $\beta^n(\Z_p\times \{0\})\cap \alpha^n(\Z_p \times\{0\})
=\Z_p \times\{0\}$,
as claimed.
\end{example}
It is unknown whether the pathology
described in Example~\ref{ex2}
can occur also if $\car(\K)=0$.
If not, we could
define a map
$\partial\hspace*{-.3mm}\Aut(G)\to \partial\hspace*{-.3mm}\Aut(L(G))$,
$\partial \alpha\mto \partial L(\alpha)$,
for each $p$-adic Lie group~$G$.
By Example~\ref{ex3} above, this map would
not always be injective.
\section{Conditions ensuring that
{\boldmath $\Phi$} is well-behaved}\label{secposi}
In this section, we describe situations
where the pathologies just encountered
can be ruled out.
In particular, we shall see that
$\Phi\colon\partial G\to\Aut(L(G))$
is a well-defined isometric
embedding for $G$ in a large
class of linear algebraic groups
(the ``generalized Cayley groups'').
Recall that
if $G$ is a totally disconnected,
locally compact group, $\alpha\in \Aut(G)$
and $V\sub G$ a compact open subgroup,
then $(\alpha^n(V))_{n\in \N_0}$ is called
the \emph{ray generated by~$\alpha$ based at~$V$}
(see \cite[p.\,394]{BaW}).
The heart of this section is a technical result
(the ``Intertwining Lemma''), which enables us to
compare rays in two Lie groups
(which need not even be locally isomorphic).\\[3mm]
As the basis for our considerations,
we need some basic facts concerning
the local structure
of Lie groups over local fields and Haar
measure on them.
The following notation will be used:
If $(E,\|.\|)$ is a normed space, $r>0$ and $x\in E$,
we write
$B_r^E(x)$ (or simply $B_r(x)$)
for the ball $\{y\in E\colon \|y-x\|<r\}$.
\begin{la}\label{localHaar}
Let $G$ be a $C^1$-Lie group
over a local field~$\K$ and $\phi\colon  P\to Q$
be a $C^1$-diffeomorphism from an open identity neighbourhood
$P\sub G$ onto an open $0$-neighbourhood~$Q$
in a finite-dimensional $\K$-vector space~$E$,
such that $\phi(1)=0$.
Let $\|.\|$ be an ultrametric norm on~$E$.
Then there exists $r>0$ such that
$B_r:=B_r^E(0)$ is contained in~$Q$, and the following holds:
\begin{itemize}
\item[\rm (a)]
$W_s:=\phi^{-1}(B_s)$ is a compact open subgroup
of~$G$, for each $s\in \;]0,r]$.
In particular, $x*y:=\phi(\phi^{-1}(x)\phi^{-1}(y))$
defines a group multiplication on~$B_r$
which makes $\phi|_{W_r}^{B_r}$
an isomorphism of $C^1$-Lie groups.
\item[\rm (b)]
$B_s$ is a normal subgroup of $(B_r,*)$,
for each $s\in \;]0,r]$.
\item[\rm (c)]
$x*B_s=B_s*x=x+B_s=B_s(x)$, for each $x\in B_r$ and $s\in \;]0,r]$.
\item[\rm (d)]
The Haar measure on $(B_r,+)$
coincides with Haar measure on~$(B_r,*)$.
\end{itemize}
\end{la}
\begin{proof}
See \cite[Proposition~2.1\,(a),\,(b)]{ANA}
for the proof of (a)--(c).\footnote{The hypotheses of
{\em loc.\,cit.}\
that $E=L(G)$ and $d\phi(1)=\id_{L(G)}$
is not used in the proof.}\vspace{2mm}

(d)\hspace*{.3mm}\footnote{Cf.\ \cite[Part~II, Chapter~IV, Exercise~5]{Ser}
for a different argument in the analytic case.}
Let $\mu$ be a Haar measure on $(B_r,+)$.
Given $x\in B_r$, we consider the left translation map
$\lambda_x\colon  B_r\to B_r$, $\lambda_x(y):=x*y$.
Then $\lambda_x(B_s(y))=x*(y*B_s)=x*y+ B_s$
and thus $\mu(\lambda_x(B_s(y)))=\mu(x*y+ B_s)=\mu(y+B_s)=\mu(B_s(y))$
for all $y\in B_r$ and $s\in \;]0,r]$.
If $U\sub B_r$ is an open subset,
then $U=\bigcup_{j\in J}B_{s_j}(y_j)$
for a countable family $(B_{s_j}(y_j))_{j\in J}$
of mutually disjoint balls
(e.g., we can take the set of all balls
$B_{s(y)}(y)$ for $y\in U$
with $s(y):=\max\{t\in \;]0,1]\colon B_t(y)\sub U\}$;
cf.\ also Theorem~1 in Appendix~2
of \cite[Part~II, Chapter~III]{Ser}).
Hence $\mu(\lambda_x(U))=\sum_{j\in J}\mu(\lambda_x(B_{s_j}(y_j)))=
\sum_{j\in J}\mu(B_{s_j}(y_j))=\mu(U)$.
Since $\mu$ is outer regular,
we deduce that $\mu(\lambda_x(A))=\mu(A)$
for each Borel set $A\sub B_r$. Hence $\mu$ is a Haar measure
on $(B_r,*)$.
\end{proof}
The next lemma is the key to our positive results.
\begin{la}[Intertwining Lemma]\label{compala}
Let $G_1$ and $G_2$ be $C^1$-Lie groups
over a local field~$\K$.
Let~$\alpha_j$ and $\beta_j$
be $C^1$-automorphisms of $G_j$
and $\Omega_j\sub G_j$
be an identity neighbourhood
such that $\alpha_j(\Omega_j)=\beta_j(\Omega_j)=\Omega_j$,\vspace{.4mm}
for $j\in \{1,2\}$.
Assume further that there exists a
map $\kappa\colon  \Omega_1\to\Omega_2$
such that $\kappa(1)=1$,
$\alpha_2\circ \kappa=\kappa\circ\alpha_1|_{\Omega_1}$,
$\beta_2\circ \kappa=\kappa\circ \beta_1|_{\Omega_1}$,
and such that $\eta:=\kappa|_R^S\colon R\to S$
is a $C^1$-diffeomorphism
for some open identity neighbourhoods
$R\sub \Omega_1$ and $S\sub \Omega_2$.
Then the following can be achieved
after shrinking~$R$ and~$S$ if necessary:
\begin{itemize}
\item[\rm (a)]
$R$ and $S$ are compact open subgroups
of $G_1$ and~$G_2$, respectively.
\item[\rm (b)]
$\eta$ takes Haar measure on $R$ to
Haar measure on~$S$.
\item[\rm (c)]
Let $\cB$ be the set of all compact open subgroups
$U$ of~$G_1$ such that $U\sub R$ and $\kappa(U)$
is a compact open subgroup of~$G_2$.
Then $\cB$ is a basis for the filter of identity neighbourhoods
in~$G_1$.
\item[\rm (d)]
Given $U\in \cB$, abbreviate $U':=\kappa(U)$. Then
$d_+(\alpha^n_1(U),\beta^k_1(U))=d_+(\alpha_2^n(U'),\beta_2^k(U'))$ and
$d_+(\beta^n_1(U),\alpha^k_1(U))=d_+(\beta_2^n(U'),\alpha_2^k(U'))$,
for all $n,k\in \Z$.
\item[\rm (e)]
If each of $\alpha_1$, $\alpha_2$, $\beta_1$ and $\beta_2$ moves to infinity,
then $\partial\alpha_1=\partial\beta_1$ if and only if
$\partial \alpha_2=\partial\beta_2$.
\item[\rm (f)]
If $s_{G_1}(\alpha_1)=s_{G_2}(\alpha_2)>1$
and $s_{G_1}(\beta_1)=s_{G_2}(\beta_2)>1$, then
$\delta_+(\alpha_1,\beta_1)
=\delta_+(\alpha_2,\beta_2)$,
$\delta_+(\beta_1,\alpha_1)=
\delta_+(\beta_2,\alpha_2)$
and $\delta(\alpha_1,\beta_1)
=\delta(\alpha_2,\beta_2)$.
\end{itemize}
\end{la}
\begin{proof}
Let $\|.\|$ be an ultrametric norm on
$E:=L(G_2)$ and
$\phi\colon  P \to Q\sub E$ be a chart
for~$G_2$ around~$1$, such that $\phi(1)=0$
and $P\sub S$.
By Lemma~\ref{localHaar}, after shrinking~$P$ and~$Q$
we may assume the following:
$Q=B_r\sub E$ for some $r>0$;
$W_s:=\phi^{-1}(B_s)$ is a compact open subgroup
of~$G_2$, for each $s\in \;]0,r]$;
and the image measure $\mu_2:=\phi^{-1}(\mu)$
is a Haar measure on~$P=W_r$,
where $\mu$ is a given Haar measure
on $(B_r,+)$.
Since also $\psi\colon  \eta^{-1}(P)\to Q$,
$\psi(x):=\phi(\eta(x))$
is a diffeomorphism
with $\psi(1)=0$, applying Lemma~\ref{localHaar}
again we see that after shrinking~$r$,
we may assume that furthermore
$V_s:=\psi^{-1}(B_s)=\eta^{-1}(W_s)$ is a compact open subgroup
of~$G_1$ for each $s\in \;]0,r]$
and $\mu_1:=\psi^{-1}(\mu)=\eta^{-1}(\mu_2)$ a Haar measure
on~$V_r$.\vspace{2mm}

(a) After replacing~$R$ and~$S$
with $V_r$ and $W_r$, respectively,
(a) holds.\vspace{2mm}

(b) By the preceding, indeed
$\eta(\mu_1)=\eta(\eta^{-1}(\mu_2))=\mu_2$.\vspace{2mm}

(c) It is clear that the sets~$V_s$
provide a basis of identity neighbourhoods,
and we have $\{V_s\colon  s\in \;]0,r]\}\sub \cB$
since $\kappa(V_s)=W_s$.\vspace{2mm}

(d) Let $U\in \cB$ and abbreviate $U':=\kappa(U)$.
Then $\kappa(\alpha^{-n}_1(\beta^k_1(U)))
=\alpha_2^{-n}(\beta_2^k (\kappa(U)))=
\alpha^{-n}_2(\beta_2^k(U'))$ for all $n,k \in \Z$,
entailing that
\begin{eqnarray*}
\hspace*{-4mm}\lefteqn{[\alpha^n_1(U):
\alpha^n_1(U)\cap \beta^k_1(U)]}\qquad\\
& = & [U: U\cap \alpha_1^{-n}\beta^k_1(U)]
\hspace*{5.5mm} =\; \frac{\mu_1(U)}{\mu_1(U\cap \alpha_1^{-n}\beta^k_1(U))}\\
& = & \!
\frac{\mu_2(\eta(U))}{\mu_2(\eta(U\cap \alpha^{-n}_1\beta^k_1(U)))}
\; =\; \frac{\mu_2(\kappa(U))}{\mu_2(\kappa(U)\cap
\kappa(\alpha^{-n}_1 \beta^k_1(U)))}\\
& = &
\frac{\mu_2(U')}{\mu_2(U'\cap \alpha_2^{-n}\beta_2^k(U'))}
\;\;\; =\;  [U': U'\cap \alpha_2^{-n}\beta_2^k(U')]\\
& = & [\alpha_2^n(U'): \alpha_2^n(U')\cap \beta_2^k(U')]
\end{eqnarray*}
and thus $d_+(\alpha_1^n(U),\beta^k_1(U))= d_+(\alpha_2^n(U'),\beta_2^k(U'))$.
Interchanging the roles of~$\alpha_j$ and~$\beta_j$,
we see that also
$d_+(\beta_1^n(U),\alpha_1^k(U))
=d_+(\beta_2^n(U'),\alpha_2^k(U'))$.\vspace{2mm}

(e) Let $U\in \cB$ and $U':=\kappa(U)$.
By (d), we have
$d(\alpha_1^n(U),\beta_1^k(U))=
d(\alpha_2^n(U'),\beta_2^k(U'))$ for all $k, n\in \N$.
Hence, if
$(d(\alpha_1^{nk}(U),\beta_1^{n\ell}(U)))_{n\in \N}$ is\linebreak
bounded for certain $k,\ell\in \N$, then so is
the sequence
$(d(\alpha_2^{nk}(U'),\beta_2^{n\ell}(U')))_{n\in \N}$,
and vice versa (as both coincide).
Thus $\partial\alpha_1=\partial\beta_1$ if and only
if $\partial\alpha_2=\partial \beta_2$.\vspace{2mm}

(f) Let $U$ and $U'$ be as before.
Then
\begin{eqnarray*}
\delta_n^{U,U}(\alpha_1,\beta_1)
\!&\!=\!&\!
\min \Big\{\frac{d_+(\alpha_1^n(U),\beta^k_1(U))}{n\log(s_{G_1}(\alpha_1))}
\colon 
\mbox{$k\in \N$ \,s.t.\
$s_{G_1}(\beta_1)^k\leq s_{G_1}(\alpha_1)^n$}\Big\}\\
\! &\!=\!&\!
\min \Big\{\frac{d_+(\alpha_2^n(U'),\beta^k_2(U'))}{n\log(s_{G_2}(\alpha_2))}
\colon 
\mbox{$k\in \N$ \,s.t.\
$s_{G_2}(\beta_2)^k\leq s_{G_2}(\alpha_2)^n$}\Big\}\\
&\!=&
\delta_n^{U',U'}(\alpha_2,\beta_2)
\end{eqnarray*}
for each $n\in \N$,
using (d) and the hypothesis
that
$s_{G_1}(\alpha_1)=s_{G_2}(\alpha_2)$
and $s_{G_1}(\beta_1)=s_{G_2}(\beta_2)$.
As a consequence, $\delta_+(\alpha_1,\beta_1)=\delta_+(\alpha_2,\beta_2)$.
The same argument
gives $\delta_+(\beta_1,\alpha_1)=\delta_+(\beta_2,\alpha_2)$,
whence also $\delta(\alpha_1,\beta_1)=\delta(\alpha_2,\beta_2)$.
\end{proof}
Given a Lie group~$G$
and $x\in G$,
we use the notation
$I_x\colon G\to G$, $y\mto xyx^{-1}$
for the inner automorphism associated with~$x$
and set $\Ad_x:=L(I_x):=T_1(I_x)$.
In the following,
we consider $G$ (and each conjugation-invariant
subset) as a $G$-space
via $x.y:=I_x(y)$.
We consider $L(G)$ as a $G$-space via
$x.y:=\Ad_x(y)$.
When speaking of a linear algebraic group~$G$
over a local field~$\K$, more precisely we
mean the Lie group of $\K$-rational points
(cf.\ \cite[Chapter~I, Proposition~2.5.2]{Mar}).
We occasionally write
$\cg:=L(G)$ for the Lie algebra
of a Lie group (or linear algebraic group)~$G$.
\begin{prop}\label{extracond}
Let $G$ be a $C^1$-Lie group
over a local field~$\K$.
Assume that there exists a map
$\kappa\colon \Omega \to L(G)$
on a conjugation-invariant
identity neighbourhood
$\Omega\sub G$
such that $\kappa(1)=0$,
$\kappa$ is $G$-equivariant
$($i.e.\ $\kappa\circ I_x|_\Omega=\Ad_x\circ \kappa$
for all $x\in G)$, and $\kappa|_R^S$
is a $C^1$-diffeomorphism for
some identity neighbourhood $R\sub \Omega$ and
some $0$-neighbourhood $S\sub L(G)$.
Then the map
\[
\Phi\colon  \partial G \to \partial \! \Aut(L(G))\, ,
\quad \Phi(\partial \alpha):=\partial L(\alpha)
\]
is well-defined and injective.
If, furthermore,
$s_G(\alpha)=s_{L(G)}(L(\alpha))$
for each $\alpha\in \Inn(G)_>$,
then $\Phi$ is an isometric embedding.
\end{prop}
\begin{rem}
If $G$ is a linear algebraic
group over~$\K$,
then $s_G(\alpha)=s_{L(G)}(L(\alpha))$
for each $\alpha\in \Inn(G)$,
as a special case of Lemma~\ref{liescale}\,(d).
\end{rem}
{\bf Proof of Proposition~\ref{extracond}.}
If $\alpha,\beta\in \Inn(G)$
move to infinity,
then also $L(\alpha)$ and $L(\beta)$
move to infinity
(see Lemma~\ref{toinfty}\,(a)).
Applying Lemma~\ref{compala}\,(e)
to the automorphisms~$\alpha$, $\beta$
and~$L(\alpha)$, $L(\beta)$ of the Lie
groups~$G$ and $(\cg,+)$, respectively,
we see that $\partial L(\alpha)=\partial L(\beta)$
if and only if $\partial \alpha=\partial \beta$.
Hence $\Phi$ is well-defined and injective.
If, furthermore,
$s_G(\alpha)=s_{L(G)}(L(\alpha))$
and $s_G(\beta)=s_{L(G)}(L(\beta))$
for all $\alpha,\beta \in \Inn(G)_>$,
then $\delta (\partial \alpha,\partial \beta)=
\delta(\partial L(\alpha), \partial L(\beta))$
by Lemma~\ref{compala}\,(f)
and thus $\Phi$ is an isometry.\,\vspace{2.3mm}\Punkt

\noindent
To illustrate Proposition~\ref{extracond},
we now consider various classes
of examples, which can be discussed by
elementary means. Afterwards,
we define a quite general class
of linear algebraic groups
with similar properties.
\begin{cor}\label{orthog}
Let $\K$ be a local field
and $n\in \N$.
Assume that $G$ is either
\begin{itemize}
\item[\rm(a)]
the general linear group $\GL_n(\K)$; or
\item[\rm(b)]
the special linear group $\SL_n(\K):=\{g\in\GL_n(\K)\colon
\dt g=1\}$, provided that
$\car(\K)=0$ or $\car(\K)>0$ and $\car(\K)$
does not divide~$n$; or
\item[\rm(c)]
the orthogonal group
$\Orth_n(\K):=\{g\in \GL_n(\K)\colon  g^T=g^{-1}\}$,
where $\car(\K)\not=2$; or
\item[\rm(d)]
the group
$\UT_n(\K):=\{(a_{ij})_{i,j=1}^n\in
\GL_n(\K)
\colon i> j\impl a_{ij}=0\}$
of all invertible upper triangular matrices.
\end{itemize}
Then $\Phi \colon \partial G \to \partial \Aut(L(G))$,
$\partial \alpha\mto \partial L(\alpha)$
is well-defined and an isometry.
\end{cor}
\begin{proof}
(a) If $G=\GL_n(\K)$, let $\cg:=\gl_n(\K)$ and
consider the map
$\kappa\colon  G \to \cg$,
$\kappa(x):=x-\one$.
Then $\kappa(\one)=0$ and
$\kappa$ is a $C^\omega$-diffeomorphism
onto the open subset $G-\one$ of~$\cg$
(with inverse $x\mto x+\one$).
Given $x,y\in G$, we have
\[
\Ad_x(\kappa(y))=x\kappa(y)x^{-1}=x(y-\one)x^{-1}
=xyx^{-1}-\one=\kappa(I_x(y))\,,
\]
whence~$\kappa$ is $G$-equivariant.
Now apply Proposition~\ref{extracond}.\vspace{2mm}

(b) The map
$\kappa\colon \SL_n(\K)\to\Sl_n(\K)$,
$\kappa(g):=g-\frac{\tr(g)}{n}\one$
is $G$-equivariant, and $\kappa(\one)=0$.
Furthermore, $d\kappa(\one)=\id\in \GL(\cg)$,
since $d\kappa(\one).\gamma'(0)=(\kappa\circ\gamma)'(0)
=\gamma'(0)+\frac{\tr\gamma'(0)}{n}\one
=\gamma'(0)$
for each analytic map $\gamma\colon U\to \SL_n(\K)$
on some $0$-neighbourhood $U\sub \K$ with $\gamma(0)=\one$.
The second summand vanishes because
$\gamma'(0)\in \Sl_n(\K)$.
Hence~$\kappa$ is a local diffeomorphism
at~$\one$, and
Proposition~\ref{extracond} applies.\vspace{2mm}

(c) If $\car(\K)\not=2$ and
$G=\Orth_n(\K)$,
we let $\cg:=L(G)=\co_n(\K)=$
$\{X\in \gl_n(\K)\colon  X^T=-X\}$
be the orthogonal Lie algebra.
Then
\begin{equation}\label{ome1}
\Omega \; := \; \{g\in \Orth_n(\K)\colon  \one+g\in \GL_n(\K)\}
\end{equation}
is an open conjugation-invariant identity neighbourhood in $G$ and
the Cayley transform
\[
\kappa\colon  \Omega \to \cg\, ,
\qquad \kappa(g)\,:=\, (\one-g)(\one+g)^{-1}
\]
is a $C^\omega$-diffeomorphism
onto an open $0$-neighbourhood in~$\cg$
(as we recall in Appendix~\ref{appA}).
We have $\kappa(\one)=0$,
and furthermore~$\kappa$ is $G$-equivariant,
since
$\Ad_x(\kappa(y))=x\kappa(y)x^{-1}=x(\one- y)(\one+y)^{-1}x^{-1}
=(\one-xyx^{-1})(\one+xyx^{-1})^{-1}=\kappa(xyx^{-1})$
for all $x\in G$ and $y\in \Omega$.
Hence Proposition~\ref{extracond}
applies.\vspace{2mm}

(d) Let $\ut_n(\K):=\{(a_{ij})_{i,j=1}^n\in M_n(\K)\colon
i> j\impl a_{ij}=0\}$ be the Lie algebra
of upper triangular matrices.
Then $\kappa\colon \UT_n(\K)\to \ut_n(\K)$,
$x\mto x-\one$ satisfies the hypotheses
of Proposition~\ref{extracond}.
\end{proof}
The following definition captures
the essence of the arguments
just used.
\begin{defn}\label{defcay}
Let $G$ be a Zariski-connected
linear algebraic group
over an infinite field~$\K$,
with Lie algebra~$\cg$.
We say that~$G$ is a
\emph{generalized}~\mbox{\emph{Cayley}} \emph{group}
if there exists a $G$-equivariant
rational map
$\kappa\colon G \dotarrow \cg$
defined over~$\K$,
such that $\kappa(1)$ is defined,\footnote{The dotted arrow
indicates that $\kappa$ is only partially
defined. For further information concerning
rational maps, cf.\ \cite{Che}.}
$\kappa(1)=0$ and $d\kappa(1)\in \GL(\cg)$.
\end{defn}
\begin{rem}\label{remcay}
Following~\cite{LPR},
$G$ is called a ($\K$-) \emph{Cayley group}
if there exists a $G$-equivariant
birational isomorphism
$\kappa\colon G\dotarrow \cg$
(defined over~$\K$).
Also some weakened versions of this concept
were considered in~\cite{LPR}.
Notably, for~$\K$ an algebraically
closed field of characteristic~$0$,
they showed that
each connected linear algebraic group~$G$
admits a $G$-equivariant and
dominant morphism $G\to \cg$
of affine varieties~\cite[Theorem~10.2]{LPR}.
Another result is more relevant for us:
For~$\K$ as before,
each connected reductive group~$G$
admits a $G$-equivariant birational isomorphism
$\kappa \colon G\to \cg$
which is a morphism of algebraic
varieties, takes~$1$ to~$0$,
and is \'{e}tale at~$1$
(see \cite[Corollary to Lemma~10.3]{LPR};
cf.\ also \cite{KaM} for related earlier
studies).
Hence, every reductive group over an algebraically closed
field of characteristic~$0$
is a generalized Cayley group
in our sense.
Unfortunately, no comparable results
seem to be available yet for algebraic groups
over ground fields
which are not algebraically
closed or have positive characteristic.
But it is to be expected
that also many of these will be generalized
Cayley groups.
The examples given in Corollary~\ref{orthog}\,(d)
show that also some
non-reductive groups are (generalized)
Cayley groups (see also \cite[Example~1.21]{LPR}).
\end{rem}
\begin{cor}\label{corcay}
If~$G$ is a generalized Cayley group
over a local field~$\K$, then
$\Phi\colon  \partial \hspace*{.4mm}G \to \partial \hspace*{-.2mm}\GL(\cg)$,
$\partial \alpha\mto \partial L(\alpha)$
is a well-defined isometric embedding.
\end{cor}
\begin{proof}
For $\kappa$ as in Definition~\ref{defcay},
its domain of definition~$\Omega$
is conjugation-invariant
and also all other hypotheses
of Proposition~\ref{extracond}
are satisfied.\vspace{-2mm}
\end{proof}
\begin{rem}
Note that the isometry $\Phi$ in Corollary~\ref{corcay}
factors to an isometry $\partial \hspace*{.3mm}G/\delta^{-1}(0)\to \partial
\hspace*{-.2mm}\GL(\cg)/\delta^{-1}(0)$ between
the corresponding metric spaces,
which in turn extends uniquely
to an isometry
$\wb{\partial \hspace*{.38mm}G} \to \wb{\partial \hspace*{-.22mm}\GL(\cg)}$
between the completions.
\end{rem}
\section{Non-existence of small joint tidy subgroups}\label{secsmall}
In this section, we provide counterexamples
showing that the existence of small tidy subgroups
for each individual automorphism in a flat group~$\cH$
need not ensure the existence of
small joint tidy subgroups for~$\cH$, not even
if $\cH$ is finitely generated.
We also show that a flat group
may contain automorphisms which are not tidy,
although it is generated by a set
of tidy automorphisms.
For general information
concerning flat groups, see~\cite{Wi2b}.\\[3mm]
Throughout this
section, $J$ is an infinite set,
$F$ a non-trivial finite group
and $G:=F^J$, equipped with the (compact) product topology.
The group
$\Sym(J)$ of all bijective self-maps
of~$J$ admits a permutation
representation $\pi\colon \Sym(J)\to \Aut(G)$
on~$G$ by automorphisms via $\pi(\sigma)(f)(j):=f(\sigma^{-1}(j))$
for $\sigma\in \Sym(J)$, $f\in G$, $j\in J$.
Furthermore, $\Sym(J)$
(and its subgroups) act on~$J$ in an obvious way.
\begin{la}\label{laprep}
Let $H\leq \Sym(J)$ be a subgroup
and $\cH:=\pi(H)\leq \Aut(G)$.
\begin{itemize}
\item[\rm(a)]
Then $G$ is tidy for each $\alpha\in\cH$,
and hence~$\cH$ is flat.
Furthermore, $s_G(\alpha)=1$ for each $\alpha\in \cH$. 
\item[\rm(b)]
A compact open subgroup $U\sub G$
is tidy for $\alpha\in \cH$ if and only if
$\alpha(U)=U$.
\item[\rm(c)]
If every $H$-orbit in~$J$ is infinite
$($\hspace{-.3mm}e.g., if $H$ acts transitively on~$J\hspace*{.4mm})$,
then $G$ is the only joint tidy subgroup
for~$\cH$.
\item[\rm(d)]
If $\sigma\in H$ and all $\langle\sigma\rangle$-orbits
in~$J$ are finite, then $\pi(\sigma)$ is
a tidy automorphism of~$G$.
\end{itemize}
\end{la}
\begin{proof}
(a) For $\alpha\in \cH$,
we have $\alpha^{-1}(G)=G$
and thus
$[G:G\cap\alpha^{-1}(G)]=1$,
which is minimal.
Hence $G$ is tidy for each $\alpha\in \cH$
and $s_G(\alpha)=1$.\vspace{1mm}

(b) Since $s_G(\alpha)=s_G(\alpha^{-1})=1$ by~(a),
it is well-known that~$U$ is tidy
if and only if $\alpha(U)=U$
holds.\footnote{This is obvious
from an alternative
definition of tidy subgroups \cite[Definition~2.1]{Wi2},
which is equivalent to the one we use
by \cite[Theorem~3.1]{Wi2}.
Cf.\ also~\cite[Corollary~3.11]{Wi2}.}\vspace{1mm}

(c) Assume that all orbits of~$H$ are
infinite and $U\sub G$ is tidy for each $\alpha\in \cH$.
Then $\{(x_j)_{j\in J}\in G\colon j\in A\impl x_j=1\}\sub U$
for some finite set $A\sub J$.
Given $x\in F$ and $k\in J$,
let $f_k(x)$ be the element
of~$G$ with $k$-th component~$x$
and all other entries~$1$.
The orbit $H.k$ being infinite, there
exists $m\in H.k \setminus A$.
Since $m\in H.k$,
there exists $\sigma\in H$
such that $\sigma(m)=k$.
Then $f_m(x)\in U$ and hence $f_k(x)=
\pi(\sigma)(f_m(x))\in \pi(\sigma)(U)=U$,
using~(b).
Hence~$U$ contains the dense
subgroup $\langle f_k(x)\colon k\in J, x\in F\rangle$ of~$G$
and thus $U=G$.\vspace{1mm}

(d) If all orbits of~$\sigma$ are finite
and $V\sub G$ is an identity neighbourhood,
pick a finite set $A\sub J$ with
$U:=\{(x_j)_{j\in J}\in G\colon j\in A\impl x_j=1\}\sub V$.
Since~$\sigma$ has finite orbits, after
increasing~$A$ we may assume that~$A$
is a union of $\sigma$-orbits. Set $\alpha:=\pi(\sigma)$.
Then $\alpha(U)=U$, and thus
$U$ is a subgroup of~$V$ which is tidy for~$\alpha$
(by~(b)).
\end{proof}
\begin{example}
In this example, we specialize to
the case $J=\Z$. Thus $G=F^\Z$.
We let $\cH:=\pi(\Sym_{\fin}(\Z))$, where
\[
\Sym_{\fin} (\Z)\; :=\;
\{\sigma \in \Sym(\Z)\colon  \mbox{$\sigma(n)=n$
for all but finitely many~$n$}\}
\]
is the group of all finite permutations of~$\Z$.
Since $\Sym_{\fin}(\Z)$ acts transitively
and each $\sigma\in \Sym_{\fin}(\Z)$
has finite orbits,
Lemma~\ref{laprep}
shows that~$\cH$ is flat,
each $\alpha\in \cH$ is tidy,
but~$G$ is the only compact open subgroup
of~$G$ which is tidy for all $\alpha\in \cH$
simultaneously.
\end{example}
There even is a finitely generated counterexample.
\begin{example}
Let $J$ be a finitely generated, algebraically periodic,
infinite group (for example,
a ``Tarski monster'' as in~\cite{Ols}).
Then~$J$ acts on itself via left translation.
We let $H\leq\Sym(J)$ be the corresponding
group of permutations and $\cH:=\pi(H)\leq \Aut(G)$.
Then $\cH$ is a finitely generated
group. Since $H$ acts transitively on~$J$
and each $\sigma\in H$ has finite orbits,
Lemma~\ref{laprep}
shows that~$\cH$ is flat and
each $\alpha\in \cH$ is tidy,
but~$G$ is the only joint tidy subgroup for~$\cH$.
\end{example}
Finally, we show that,
if a flat group~$\cH$ of automorphisms
is generated by a set of tidy automorphisms,
this does not imply that
each $\alpha\in \cH$ is tidy.
\begin{example}
Take $J:=\Z$;
thus $G=F^\Z$.
We define $\sigma, \tau\in \Sym(\Z)$
via $\sigma(2k+1):=2k$,
$\sigma(2k):=2k+1$,
$\tau(2k-1):=2k$ and $\tau(2k):=2k-1$
for $k\in \Z$.
We also set $H:=\langle \sigma,\tau\rangle\leq\Sym(\Z)$
and $\cH:=\pi(H)$.
Then~$\cH$ is flat, by Lemma~\ref{laprep}\,(a).
Furthermore,
$\pi(\sigma)$ and $\sigma(\tau)$ are tidy,
since all orbits of~$\sigma$ and~$\tau$ have
two elements
(see Lemma~\ref{laprep}\,(d)).
However, although both generators~$\pi(\sigma)$ and~$\pi(\tau)$
of~$\cH$
are tidy,
the group~$\cH$ contains automorphisms
which are not tidy.
For example, consider the element
$\alpha:=\pi(\sigma\circ\tau)\in \cH$.
Since $(\sigma\circ \tau)(2k)=2k-2$
and $(\sigma\circ\tau)(2k-1)=2k+1$
for each $k\in \Z$,
we see that
$2\Z$ and $2\Z+1$
are the orbits of $\sigma\circ\tau$.
Since both of these are infinite,
the only subgroup of~$G$ tidy
for each $\beta\in \langle\alpha \rangle$
is all of~$G$, by Lemma~\ref{laprep}\,(c).
Since a subgroup is tidy for~$\alpha$
if and only if it is tidy for each
$\beta\in \langle\alpha\rangle$
(cf.\ Lemma~\ref{laprep}\,(b)),
we deduce that~$G$ is the only subgroup of~$G$
tidy for~$\alpha$, and thus~$\alpha$ is not tidy.
\end{example}
\appendix
\section{Basic properties of the Cayley transform}\label{appA}
For the convenience of the reader,
we provide proofs for the following (well-known) properties of the Cayley
transform which we used in the proof
of Corollary~\ref{orthog}\,(c)
(cf.\ \cite[pp.\,123--125]{Che} for further information).
\begin{la}\label{Ct}
Given $n\in \N$ and a local field $\K$ such that $\car(\K)\not=2$,
define $\Omega:=\{x\in M_n(\K)\colon  \one+x\in \GL_n(\K)\}$.
Then~$\Omega$ is an open subset of $M_n(\K)$ such that $\{{\bf 0},\one\}
\sub \Omega$, and
\[
\theta\colon  \Omega\to\Omega\, ,
\quad \theta(x)\, :=\, (\one-x)(\one+x)^{-1}
\]
is a bijective self-map of~$\Omega$
such that $\theta\circ\theta=\id_\Omega$.
The following holds:
\begin{itemize}
\item[\rm (a)]
$\one+\theta(x)\in \GL_n(\K)$
for each $x\in \Omega$, with
$(\one+\theta(x))^{-1}=\frac{1}{2}(\one+x)$.
\item[\rm (b)]
$\Omega_1 := \Omega\cap \Orth_n(\K)$
is an open conjugation-invariant identity neighbourhood
in~$\Orth_n(\K)$ and
$\Omega_2 := \Omega\cap \co_n(\K)$
is an open $0$-neighbourhood in~$\co_n(\K)$,
such that $\theta(\Omega_1)=\Omega_2$.
\item[\rm(c)]
The map $\theta|_{\Omega_1}^{\Omega_2}\colon \Omega_1\to\Omega_2$
is a $C^\omega$-diffeomorphism.
\end{itemize}
\end{la}
\begin{proof}
(a) We have $(\one+\theta(x))(\one+x)=(\one+(\one-x)(\one+x)^{-1})(\one+x)
=\one+x+\one-x=2\one$ for each $x\in \Omega$,
showing that indeed $\one+\theta(x)\in \GL_n(\K)$
(and thus $\theta(x)\in \Omega$), with
$(\one+\theta(x))^{-1}=\frac{1}{2}(\one+x)$.\vspace{2mm}

{\em $\theta$ is an involution\/} as
$\theta(\theta(x))=(\one-(\one-x)(\one+x)^{-1})(\one+\theta(x))^{-1}
=\frac{1}{2}(\one-(\one-x)(\one+x)^{-1})(\one+x)
=\frac{1}{2}(\one+x-(\one-x))=x$ for each $x\in \Omega$, using\,(a).
In particular, this implies
that $\theta\colon \Omega\to\Omega$ is a bijection.\vspace{2mm}

(b)
If $g\in \Omega_1$,
then $\theta(g)^T=(\one+g^T)^{-1}(\one-g^T)=
(\one+g^{-1})^{-1}(\one-g^{-1})
=(g^{-1}(g+\one))^{-1}g^{-1}(g-\one)
=(g+\one)^{-1}(g-\one)=-\theta(g)$.
Thus $\theta(g)\in \co_n(\K)\cap\Omega=\Omega_2$.
Conversely, $\theta(x)^T=
(\one+x^T)^{-1}(\one-x^T)=(\one- x)^{-1}(\one+x)$
for $x\in \Omega_2$
and thus $\theta(x)^T\theta(x)=
(\one- x)^{-1}(\one+x)(\one-x)(\one+x)^{-1}=\one$,
using that all matrices involved commute.
Thus $\theta(x)\in \Omega\cap \Orth_n(\K)=\Omega_1$.
Now~(b) follows.\vspace{2mm}

(c) Since $\theta$ is $C^\omega$,
also $\theta|_{\Omega_1}^{\Omega_2}$
and $\theta|_{\Omega_2}^{\Omega_1}=(\theta|_{\Omega_1}^{\Omega_2})^{-1}$
are~$C^\omega$.
\end{proof}
{\footnotesize
{\bf Helge Gl\"{o}ckner}, TU Darmstadt, Fachbereich Mathematik AG~5,
Schlossgartenstr.\,7,\\
64289 Darmstadt, Germany.\\
E-Mail:
\,{\tt gloeckner@mathematik.tu-darmstadt.de}\\[3mm]
{\bf George A. Willis}, School of Mathematical and Physical Sciences,
University of\\
Newcastle, University Drive, Callaghan, NSW 2308,
Australia.\\
E-Mail: \,{\tt George.Willis@newcastle.edu.au}}
\end{document}